\documentclass[12pt, reqno]{amsart}
\usepackage{amsmath, amsthm, amscd, amsfonts, amssymb, graphicx, color}
\usepackage[bookmarksnumbered, colorlinks, plainpages]{hyperref}

\newtheorem{Theorem}{Theorem}[section]
\newtheorem{Definition}[Theorem]{Definition}
\newtheorem{Proposition}[Theorem]{Proposition}

\newtheorem{Lemma}[Theorem]{Lemma}
\newtheorem{Corollary}[Theorem]{Corollary}
\theoremstyle{remark}
\newtheorem{Example}[Theorem]{Example}
\newtheorem{Remark}[Theorem]{Remark}

\def\ovr{\overline}

\def\om{\omega}
\def\Om{\Omega}

\def\sbs{\subset}

\def\sgn{\operatorname{sgn}}

\def\lim{\operatorname{lim}}

\def\be{\begin{enumerate}}
\def\ee{\end{enumerate}}
\def\bT{\begin{Theorem}}
\def\eT{\end{Theorem}}
\def\bP{\begin{Proposition}}
\def\eP{\end{Proposition}}
\def\bPr{\begin{proof}}
\def\ePr{\end{proof}}
\def\bD{\begin{Definition}}
\def\eD{\end{Definition}}
\def\bE{\begin{Example}}
\def\eE{\end{Example}}
\def\bL{\begin{Lemma}}
\def\eL{\end{Lemma}}
\def\bC{\begin{Corollary}}
\def\eC{\end{Corollary}}

\def\E{{\mathcal E}}

\def\M{{\mathcal M}}

\begin{document}
\setcounter{page}{1}

\title[Weighted Hardy spaces]{Structure of weighted Hardy spaces in the plane}

\author[N\.{i}hat G\"{o}khan G\"{o}\u{g}\"{u}\c{s}]{N\.{i}hat G\"{o}khan G\"{o}\u{g}\"{u}\c{s}$^1$ $^{*}$}

\address{$^{1}$ Sabanci University,  Tuzla , Istanbul 34956 Turkey.}
\email{\textcolor[rgb]{0.00,0.00,0.84}{nggogus@sabanciuniv.edu}}


\subjclass[2010]{Primary 30H10, 30J99, Secondary 46E22.}

\keywords{Weighted Hardy spaces, subharmonic exhaustion, dual
space.}


\date{Received: xxxxxx; Revised: yyyyyy; Accepted: zzzzzz.
\newline \indent $^{*}$ Corresponding author}

\begin{abstract}
We characterize certain weighted Hardy spaces on the unit
disk and completely describe their dual spaces.
\end{abstract} \maketitle

\section{Introduction and preliminaries}
\par By a recent paper of Poletsky and Stessin \cite{PolSte08} to each subharmonic function on a
bounded regular domain $G$ which is continuous near the boundary
corresponds a space $H^p_u$ of analytic functions in $G$ with a
certain growth condition. These are namely Poletsky-Stessin Hardy
spaces. They include and generalize the well-known classical Hardy
spaces. This new theory unifies the standpoints of various analytic
function spaces into one.

\par The first generalizations in this direction of the theory of Hardy
spaces on hyperconvex domains in ${\mathbb{C}}^n$ was suggested and
studied in \cite{Al03}. More recently the theory is extended to
hyperconvex domains in \cite{PolSte08}. Boundedness and compactness
of the composition operators on these new Poletsky-Stessin Hardy and
Bergman type spaces were investigated there. After this motivating
work more investigation \cite{AlGog}, \cite{Shr}, \cite{Sahin}
revealed the structure and first examples of these Hardy type spaces
in the plane.

\par In \cite{AlGog} to understand the scale of weighted Hardy spaces $u\to H^p_u$
Alan and the author completely characterized $H^p_u$ spaces in the
plane domains by their boundary values or by possessing a harmonic
majorant with a certain growth (see also \cite{Shr}, \cite{Sahin}).
Basically the version of the Beurling's theorem proved in
\cite{AlGog} states that to each subharmonic exhaustion $G$
corresponds an outer function $\varphi$ which belongs to the class
$H^p_u$ so that $H^p_u$ isometrically equals to $\M_{\varphi ,p}$
for $p>0$, where $\M_{\varphi ,p}$ is the space $\varphi^{2/p}H^p$
endowed with the norm
\[\|f\|_{\M_{\varphi ,p}}:=\|f/\varphi^{2/p}\|_p, \,\,\,\,f\in \M_{\varphi ,p}.\]
This result is especially useful to construct examples of analytic
function spaces enjoying certain desired properties. The space
$\M_{\varphi ,2}$, when $\|\varphi\|_{\infty}\leq 1$, was studied as
a tool to understand certain sub-Hardy Hilbert spaces in the unit
disk in \cite{Sarason}. Two problems were not answered in
\cite{AlGog}: \be
\item Can we go back? That is, given analytic $\varphi$ can one
find a subharmonic exhaustion $u$ so that $H^p_u=\M_{\varphi ,p}$?
\item  For the space $H^p_u$, consider the class of all representatives, i.e., subharmonic exhaustions $v$
so that $H^p_u=H^p_v$. What kind of "good" representatives are
there?\ee

\par In this note we give answers for both questions. We show
under certain growth conditions on the analytic function $\varphi$
on the disk that it is possible to construct a subharmonic
exhaustion $u$ on the disk so that $H^p_u$ equals to $\M_{\varphi
,p}$. Moreover, by the
construction, $u$ is real analytic and satisfies the bi-Laplacian in
the unit disk. This is a new information related to the second
question.

\par In addition, using the boundary value characterization from \cite{AlGog} we completely characterize the
dual space of $H^p_u$ and discuss the corresponding extremal and
dual extremal problems.

\par After several months of submission of this paper there appeared a preprint \cite{PS}. Theorem 3.3 in this paper
is similar to Theorem 2.1 below. In Theorem 2.1 we do not require any integrability condition on the subharmonic exhaustion
$u$, however the authors in \cite{PS} require $u$ to be intagrable.


\par Let us start to recall basic definitions. A function $u\leq 0$ on a bounded open set $G\subset \mathbb{C}$
is called an exhaustion on $G$ if the set
\[
B_{c,u}:=\{z\in G:u(z)< c\}
\]
is relatively compact in $G$ for any $c<0$. When $u$ is an
exhaustion and $c<0$, we set
\[
u_c:=\max\{u,c\},\,\,\,\,\,\, S_{c,u}:=\{z\in G:u(z)=c\}.
\]
Let $u \in sh(G)$ be an exhaustion function which is continuous with
values in $\mathbb{R}\cup \{-\infty\}$. Following Demailly
\cite{Dem87} we define
\begin{eqnarray*}
 \mu_{c,u}:=\Delta u_c-\chi_{G\backslash B_{c,u}}\Delta u,
\end{eqnarray*}
where $\chi_{\om}$ is the characteristic function of a set $\om\sbs
G$. We denote the class of negative subharmonic exhaustion functions
on $G$ by $\mathcal{E}(G)$. The class of all functions
$u\in\mathcal{E}(G)$ for which $\int\Delta u<\infty$ is denoted by
$\mathcal{E}_0(G)$.

If $u\in\mathcal{E}(G)$, then the Demailly-Lelong-Jensen formula
(\cite{Dem87}) takes the form
\begin{equation}
\int_{S_{c,u}}v\,d\mu _{c,u}=\int_{B_{c,u}}(v\Delta u-u\Delta
v)+c\int_{B_{c,u}}\Delta v,  \label{Eq:DemJenLel}
\end{equation}%
where $\mu _{c,u}$ is the Demailly measure which is supported in the level sets $%
S_{c,u}$ of $u$ and $v\in $ $sh(G)$. Let us recall that by
\cite{Dem87} if $\int_{G}\Delta u<\infty $, then the
measures $\mu _{c,u}$ converge as $c\rightarrow 0$ weak-$\ast $ in $C^*(\ovr{G})$ to a measure $%
\mu _{u}$ supported in the boundary $\partial G$.

Following \cite{PolSte08} we set
\[
\mathnormal{sh}_u(G):=\mathnormal{sh}_u:=\left\{v\in sh(G):v\geq
0,\,\sup _{c<0}\int_{S_{c,u}}v\,d\mu _{c,u}<\infty\right\},
\]%
and
\[
H^p_u(G):=H^p_u:=\left\{f\in hol(G):|f|^p\in
\mathnormal{sh}_u\right\}
\]
for every $p>0$. We write
\begin{eqnarray}\label{Eq:NormSh}
\|v\|_{u}:=\sup _{c<0}\int_{S_{c,u}}v\,d\mu _{c,u}=\int_{G}(v\Delta
u-u\Delta v)
\end{eqnarray}
for a nonnegative function $v\in sh (G)$ and set
\begin{eqnarray}
\|f\|_{u,p}:=\sup_{c<0}\left (\int_{S_{c,u}}|f|^p\,d\mu _{c,u}\right
)^{1/p}
\end{eqnarray}
for a holomorphic function $f$ on $G$. We will use
$\|f\|_u=\|f\|_{u,p}$ when $p=1$. By Theorem 4.1 of \cite{PolSte08},
$H^p_u$ is a Banach space when $p\geq 1$.  It is clear that the
function $f\equiv 1$ belongs to $H^p_u$ if and only if the Demailly measure
$\mu_{u}$ has finite mass. If $G$ is a regular bounded domain in
$\mathbb{C}$ and $w\in G$, then the Green function $v(z)=g_G(z,w)$
is a subharmonic exhaustion function for $G$. For example, when $G$
is the unit disk and $v(z)=\log |z|$, then $\mu_v$ is the normalized
arclength measure on the unit circle. We denote by $P_G(z,w)$ the
Poisson kernel for the domain $G$.

The following Theorems are recollections from \cite{AlGog}.

\bT \cite[Theorem 2.3]{AlGog}
\label{T:CharacHardyClassesHarmonicMajorant} Let $G$ be a bounded
domain, $v\geq 0$ be a function on $G$, $p>0$, and
$u\in\mathcal{E}(G)$. The following statements are equivalent:
\begin{itemize}
\item[i.] $v\in sh_u(G)$.
\item[ii.] The least harmonic majorant $h=P_G(v)$ of $\varphi$ in $G$
belongs to the class $sh_u$.
\end{itemize}
Furthermore,
\begin{eqnarray*}
\|v\|_{u}=\int _G h\Delta u=\|h\|_u.
\end{eqnarray*}
\eT

We will denote by $H^p(G)$ the space of analytic functions $f$ in
$G$ for which $|f|^p$ has a harmonic majorant in $G$ (see for
example \cite{Duren}). We always have $H^p_u\sbs H^p$ by
\cite{PolSte08}. We will denote by $\nu$ the usual arclength measure
on $\partial G$ normalized so that $\nu(\partial G)=1$. \bT
\cite[Theorem 2.10]{AlGog} \label{T:CharacHardyClassesBoundaryVal}
Let $G$ be a Jordan domain with rectifiable boundary or a bounded
domain with $C^2$ boundary, $p>1$, and $u\in\mathcal{E}(G)$. The
following statements are equivalent:
\begin{itemize}
\item[i.] $f\in H^p_u(G)$.
\item[ii.] $f\in H^p(G)$ and
$|f^*|\in L^p(V_u\nu)$,  where
\begin{eqnarray}\label{Eq:DefnV} V_u(\zeta):=\int_{G}P_G(z,\zeta)\Delta u(z),\,\,\,\,
\zeta\in\partial G.\end{eqnarray}
\item[iii.] $f\in H^p(G)$ and
there exists a positive measure $\widetilde{\mu_u}$ on $\partial G$
such that $|f^*|\in L^p(\widetilde{\mu_u})$. Moreover, if $E$ is any
Borel subset of $\partial G$ with measure $\nu(E)=0$, then
$\widetilde{\mu_u}(E)=0$ and we have the equality
\begin{eqnarray}\label{Eq:OnV} \int_{\partial G}\gamma\,d\widetilde{\mu_u}=\int_{G}P_G(\gamma)\Delta u\end{eqnarray}
for every $\gamma\in L^1(\nu)$.
\end{itemize}
In addition, if $f\in H^p_u(G)$, then
$\|f\|_{u,p}=\|f^{\ast}\|_{L^p(\widetilde{\mu_u})}$ and
$d\widetilde{\mu_u}=V_ud\nu$. \eT

\begin{Remark}\label{Rem:NormalDer} i. Theorem \ref{T:CharacHardyClassesBoundaryVal} is
valid when $p>0$ and $G$ is the unit disk or more generally a Jordan
domain with rectifiable boundary. In this case the Poisson integral
of an $L^p(d\nu)$ function $u$ has non-tangential limits equal to
$u$ on $\nu$-almost every boundary point. This is indeed what is
needed in the proof of Theorem
\ref{T:CharacHardyClassesBoundaryVal}.

ii. By replacing the function $u$ by a suitable positive multiple
$tu$, $t>0$, we may assume that $V_u\geq 1$ on $\partial G$. To do
this it is enough to take a compact set $K\sbs G$ so that $\Delta
u(K)>r>0$. Let $m:=\min_{\zeta\in\partial G}\min _{z\in
K}P_G(z,\zeta)$. Then take $t:=1/(rm)$. We will use the assumption
that $V_u\geq 1$ when convenient.

iii. The weight function $V_u$ is lower semicontinuous. To see this,
suppose $\zeta_j\in\partial G$, $\zeta_j\to\zeta$. By Fatou's lemma
\[\liminf_j V_u(\zeta_j)=\liminf_j\int_{G}P_G(z,\zeta_j)\Delta u(z)\geq \int_{G}P_G(z,\zeta)\Delta u(z)=
V_u(\zeta).\] Note that $V_u$ is the \emph{balayage} of the measure
$\Delta u$ on $\partial G$.

iv. Suppose $G$ is a bounded domain with $C^2$ boundary or a Jordan
domain with rectifiable boundary and $u\in\E_0 (G)$. Then
\[u(z)=\int_Gg_G(z,w)\Delta u(w),\,\,\,\,z\in G.\] Since
\[\frac{\partial g_G(\zeta,w)}{\partial n}=P_G(w,\zeta)\] when $\zeta\in\partial G$ and $w\in G$, $\frac{\partial
u}{\partial n}(\zeta)$ exists for every $\zeta\in\partial G$,
where $\frac{\partial }{\partial n}$ denotes the normal derivative
in the outward direction on $\partial G$ and
\[\frac{\partial u (\zeta)}{\partial n}=V_u(\zeta)=\int_{G}P_G(w,\zeta)\Delta u(w),\,\,\,\,
\zeta\in\partial G.\]

By property (\ref{Eq:OnV}) in Theorem
\ref{T:CharacHardyClassesHarmonicMajorant}
\[\int_{\partial G}V_u(\zeta)d\nu (\zeta)=\int_{G}\Delta u=\int_{\partial G}\frac{\partial
u}{\partial n}(\zeta)d\nu (\zeta).\]
\end{Remark}

\par To obtain Fatou's type results we would like to compute the
Radon-Nikodym derivative of the Demailly measures with respect to
the usual arclength measure on the level sets. In the next result we
provide this. Let $\nu_c$ denote the arclength measure on $S_{c,u}$.
Define
\begin{eqnarray*}V_{c,u}(\zeta):=\int_{B_{c,u}}P_{B_{c,u}}(z,\zeta)\Delta
u(z),\,\,\,\,\zeta\in S_{c,u},
\end{eqnarray*} where $P_{B_{c,u}}(z,\zeta)$ denotes the Poisson kernel for
$B_{c,u}$.

\bP\label{Prop:AbsContLebesgueMeas} Let $u\in\mathcal{E}(G)$, where
$G$ is a bounded regular domain. Suppose that $u$ is Lipschitz in
every compact subset of $G$. Then the measures $\nu_c$ and
$\mu_{c,u}$ are mutually absolutely continuous and $\mu
_{c,u}=V_{c,u}\nu _{c}$ with $V_{c,u}\in L^1(\nu_c)$. Moreover, for
each $c<0$ there is a constant $k_c>0$ so that $V_{c,u}\geq k_c$ on
$S_{c,u}$.\eP

\begin{proof} Let $\varphi$ be a continuous function on $S_{c,u}$ and let $h(z)$ be the harmonic function in $B_{c,u}$ with boundary
values equal to $\varphi$. By equality (\ref{Eq:DemJenLel}) we have
\begin{eqnarray*}
\int_{S_{c,u}}\varphi (\zeta )d\mu _{c,u}(\zeta )
&=&\int_{B_{c,u}}h(z)\Delta u(z) \\
&=&\int_{S_{c,u}}\left( \int_{B_{c,u}}P_{B_{c,u}}(z,\zeta )\Delta
u(z)\right)
\varphi (\zeta )\,d\nu _{c}(\zeta ) \\
&=&\int_{S_{c,u}}\varphi (\zeta )V_{c,u}(\zeta )\,d\nu _{c}(\zeta ).
\end{eqnarray*}
Hence $\mu _{c,u}=V_{c,u}\nu _{c}$. Another observation using
Fubini's theorem gives
\begin{eqnarray*}
\int_{S_{c,u}}V_{c,u} (\zeta )d\nu _{c}(\zeta ) =\int_{B_{c,u}}
\Delta u(z) =\|\mu_{c,u}\|<\infty.
\end{eqnarray*}
Thus $V_{c,u}\in L^1(\nu_c)$. Note that $\nu_c\leq k'_c\mu_{c,u}$
for some positive constant $k'_c$ by \cite{Dem87}. Hence
$V_{c,u}\geq k_c$ on $S_{c,u}$ for some $k_c>0$. This completes the
proof.
\end{proof}

\begin{Remark}
The requirement that $u$ is Lipschitz is only needed to write the harmonic measure on $B_{c,u}$
of the form $P_{B_{c,u}}d\nu_c$. There are much weaker conditions on domains for which the harmonic measure is absolutely continuous.
\end{Remark}

The next auxiliary result allows one to compare the Demailly
measures on $S_{c,u}$ with a measure on an arbitrary level set.

\bP\label{Prop:SmoothLevelSets} Let $u$ be a subharmonic exhaustion
function on a bounded regular domain $G$ in $\mathbb{C}$. Let $G_j$
be relatively compact regular open sets in $G$ so that
$\ovr{G_j}\sbs G_{j+1}$ and $\cup G_j=G$. Then for each $j$ there is
a $u_j\in\mathcal{E}(G_j)$ and for each $c<0$ there is a number $s$
with $c<s<0$ so that for any nonnegative function $v\in sh(G)$, the
integrals $\mu_{u_j}(v)$ are increasing and
\[\mu_{c,u}(v)\leq \mu_{u_j}(v)=\|v\|_{u_j}\leq
\mu_{s,u}(v).\] This means $\|v\|_u=\lim_j\mu_{u_j}(v)$ for every
nonnegative subharmonic function $v$ on $G$. \eP

\begin{proof}
Set $u_j:=u-P_{G_j}u$. Clearly $u_j\in\mathcal{E}(G_j)$. Take an
integer $j_0\geq 1$ and a number $s<0$ with $c<s$ so that
$B_{c,u}\sbs G_{j_0}\sbs B_{s,u}$. The comparison follows from
(\ref{Eq:DemJenLel}) and (\ref{Eq:NormSh}) if we note that $c\leq
P_{G_j}u$ on $B_{c,u}$ and $P_{G_j}u\leq s$ on $B_{s,u}$.
\end{proof}

If $\varphi$ is a nonzero analytic function on $\mathbb{D}$, let
$\M_{\varphi ,p}$ denote the space $\varphi^{2/p}H^p$ endowed with
the norm
\[\|f\|_{\M_{\varphi ,p}}:=\|f/\varphi^{2/p}\|_p, \,\,\,\,f\in \M_{\varphi ,p}.\]
We will call a function $\varphi\in H^2_u$ a $u$-inner function if
$|\varphi^{\ast} (\zeta)|^2V_u(\zeta)$ equals $1$ for almost every
$\zeta\in\partial\mathbb{D}$. If, moreover, $\varphi (z)$ is
zero-free, we will say that $\varphi$ is a singular $u$-inner
function. The next result is Theorem 3.2 and Corollary 3.3 from
\cite{AlGog}.

\bT\label{Th:BeurlingThm} Let $Y\not =\{0\}$ be a closed
$M_z$-invariant subspace of $H^2_u(\mathbb{D})$. Then there exists a
function $\varphi \in H^2_u$ so that $|\varphi^{\ast}
(\zeta)|^2V_u(\zeta)=1$ for almost every
$\zeta\in\partial\mathbb{D}$ and $Y=\M_{\varphi ,2}$. In particular,
there exists a $u$-inner and an outer function $\varphi \in H^2_u$
so that $H^2_u=\M_{\varphi ,2}$ and these spaces are isometric.\eT

This function $\varphi$ is determined uniquely up to a unit
constant. Note that
\begin{eqnarray}\label{Eq:VSign} V_u(e^{i\theta})=\frac{1}{|\varphi
(e^{i\theta})|^2}=\frac{1}{\varphi ^2 (e^{i\theta})}\sgn
\frac{1}{\varphi ^2 (e^{i\theta})},\end{eqnarray} where we set $sgn
\alpha:=|\alpha|/\alpha$ for any complex number $\alpha\not =0$ and
$sgn 0:=0$. If $V\geq 1$ on $\partial\mathbb{D}$, then
$|\varphi(\zeta)|\leq 1$ for almost every $\zeta$. Suppose that
$\int\Delta u<\infty$. Then the function $1$ belongs to $H^p_u$.
Hence $\varphi^{-1}$ belongs to $H^2$. Then it is an easy exercise
to show that $\varphi$ is an outer function.

\bT \label{Th:IsometryLpuLp}The set $L^p(V_ud\theta)$ coincides with
$\varphi^{2/p}L^p(d\theta)$ and the map $f\mapsto \varphi^{-2/p}f$
is an isometric isomorphism from the space $L^p(V_ud\theta)$ onto
$L^p(d\theta)$. \eT

\bT\cite[Theorem 3.4]{AlGog}\label{T:InnerOuterFactorization}
Suppose $0<p<\infty$, $f\in H^p_u(\mathbb{D})$, $f\not \equiv 0$,
and $B$ is the Blaschke product formed with the zeros of $f$. Then
there are zero-free $\varphi\in H^2_u\cap H^{\infty}$, $S\in H^{\infty}$ and $F\in H^p$ so
that $\varphi$ is outer and singular $u$-inner, $S$ is singular
inner, $F$ is outer, and
\begin{eqnarray}
 f=BS\varphi ^{2/p}F.
\end{eqnarray}
Moreover, $\|f\|_{p,u}=\|F\|_p$ and $H^p_u(\mathbb{D})=\M_{\varphi ,p}$. \eT

\bC \label{Cor:IsometryHpuHp}The map $f\mapsto \varphi^{-2/p}f$ is
an isometric isomorphism from the space $H^p_u$ onto $H^p$. \eC

The following Lemma will be useful in the next section. Its proof is a simple calculation and we outline it here.
\bL\label{Lem:DefnKappa}
Let $c$ be a number with $-1<c<0$. Then there exists a function $\kappa=\kappa_c$ defined on
$(-\infty,0]$ with the following properties:
\begin{itemize}
\item[i.] $\kappa:(-\infty,0]\to(-\infty,0]$ is non-decreasing, convex and $C^\infty$,
\item[ii.] $\kappa$ is real-analytic in $(c,0]$,
\item[iii.] $\kappa (t)\equiv c$ when $t\leq c$, $\kappa (0)=0$, and $\kappa '(0)=1$.
\end{itemize}
\eL
\begin{proof}
Let $a:=-\frac{\ln (-c)}{e}$, $b:=\frac{-1}{\ln (-c)}$, and
\[
\kappa (t):=
  \begin{cases}
      \hfill c+e^{\frac{-a}{(t-c)^b}},    \hfill & \text{ $t>c$, } \\
      \hfill c, \hfill & \text{ $t\leq c$. } \\
  \end{cases}
\]
Then
\[
\kappa '(t)=\frac{1}{e(t-c)^{b+1}}e^{\frac{-a}{(t-c)^b}}
\]
and
\[
\kappa ''(t)=\frac{1}{e(t-c)^{2b+2}}(1/e-(b+1)(t-c)^{b+1})e^{\frac{-a}{(t-c)^b}}
\]
for $t>c$. For $t\leq c$, $\kappa '(t)=\kappa ''(t)=0$. It can be checked that $\kappa ''(t)>0$
for $c<t\leq 0$, and $\kappa$ satisfies all properties in i., ii. and iii.
\end{proof}


\section{Finding subharmonic exhaustion}
Theorem \ref{T:CharacHardyClassesBoundaryVal} describes the weight
function $V_u$ corresponding to the Hardy space $H^p_u$ when the
Laplacian of $u$ is known. In Theorem
\ref{T:InnerOuterFactorization} we obtain a canonical factorization
for functions in $H^p_u$ and we see that this space is a certain
multiple $\varphi ^{2/p}H^p$ of $H^p$. The singular $u$-inner
function $\varphi$ appearing in this factorization is related to the
weight $V_u$ by
\begin{eqnarray}\label{Eq:RelationVPhi}
V_u(e^{i{\theta}})=\frac{1}{|\varphi(e^{i{\theta}})|^2},\,\,\,\,
a.e.\,\,\theta.
\end{eqnarray}
In this section we seek a converse to these results.

\par Let $G$ be a Jordan domain with rectifiable boundary and $\psi$
be a given analytic function in $H^1(G)$. The problem is to find a
subharmonic exhaustion $u$ on $G$ so that $V_u(\zeta)=|\psi
(\zeta)|$ when $\zeta\in \partial G$. Taking a conformal map of $G$
onto $\mathbb{D}$ we can always suppose that $G=\mathbb{D}$. This is
a type of \emph{inverse balayage} problem. We solve this next.
\bT
Let $\psi$ be a lower semicontinuous function on $\partial\mathbb{D}$ so that
$\psi\geq c$ for some constant $c>0$. Then there exists a function $u\in\E$ so that
$\psi=V_u$. Moreover we have the following properties:
\begin{itemize}
\item[a.] $u$ is the decreasing limit of functions in $\E_0\cap C^{\infty}(\ovr{\mathbb{D}})$ converging uniformly to $u$ on $\ovr{\mathbb{D}}$.
\item[b.] $u\in\E_0(\mathbb{D})$ if and only if $\psi\in L^1(d\nu)$.
\item[c.] If $\psi$ is $C^k$, $0\leq k\leq\infty$, on $\partial\mathbb{D}$, then $u$ is  $C^k$ on $\ovr{\mathbb{D}}$. If $\psi$ is real-analytic, then there exists a compact $K$ so that $u$ is real-analytic on $\ovr{\mathbb{D}}\backslash K$.
\end{itemize}
\eT
\begin{proof}
Suppose first that $\psi$ is $C^2$ on $\partial\mathbb{D}$ and let $\rho (re^{i\theta}):=\frac{1}{2}(r^2-1)\psi(e^{i\theta})$ for $re^{i\theta}\in\mathbb{D}$. Computing the Laplacian of $\rho$ we get
\[\Delta \rho(re^{i\theta})=2\psi (e^{i\theta})+\frac{r^2-1}{2r^2}\frac{d^2\psi(e^{i\theta})}{d\theta^2}.\]
By assumption $\Delta \rho(e^{i\theta})=2\psi (e^{i\theta})\geq 2c>0$. Hence there exists a compact $B\sbs\mathbb{D}$ so that $\Delta \rho(z)>0$ on the open set $\Om:=\mathbb{D}\backslash B$. Hence $\rho$ is a non-positive subharmonşc function on
$\Om$ and $\rho|_{\partial\mathbb{D}}\equiv 0$. Since $\rho$ is continuous on $\ovr{\mathbb{D}}$, there exists a constant $c<0$ so that the set $B_{c,\rho}$ is relatively compact in $\mathbb{D}$ and $S_{c,\rho}\sbs \Om$. Let $\kappa=\kappa_c$ be the function proivided in Lemma \ref{Lem:DefnKappa}.
Define $u(z):=\kappa (\rho(z))$ for $z\in\mathbb{D}$. Now $u=\kappa(\rho)$ is subharmonic in $\Om$, $u\equiv c$ on $B_{c,\rho}$ and $u\geq c$ on $\mathbb{D}\backslash B_{c,\rho}\sbs\Om$. Hence $u\in\E$ and $V_u=\frac{\partial u}{\partial r}=\kappa'(0)\frac{\partial \rho}{\partial r}=\psi$ on $\partial\mathbb{D}$.

Now let $\psi$ be lower semicontinuous. There exists $\psi_n$, all $C^{\infty}$ on $\partial \mathbb{D}$ so that $c\leq\psi_n(\zeta)\leq \psi_{n+1}(\zeta)$, and $\psi(\zeta)=\lim_n\psi_n(\zeta)$ for every
$\zeta\in\partial\mathbb{D}$. We let $\psi_0\equiv 0$. Replacing $\psi_n$ by $\psi_n-2^{-n}$ we may assume that $d_n:=\psi_{n+1}-\psi_n\geq 2^{-n-1}$. As in the first part of the proof we let $\rho_n(z):=\frac{1}{2}(r^2-1)d_n(e^{i\theta})$. There exists a compact $B_n\sbs\mathbb{D}$ so that $\Delta \rho_n(z)>0$ on the open set $\Om_n:=\mathbb{D}\backslash B_n$. This time we choose constants $-2^{-n}\leq c_n<0$ so that $B_{c_n,\rho_n}$ is relatively compact in $\mathbb{D}$, $S_{c_n,\rho_n}\sbs \Om_n$, and $B_{c_n,\rho_n}\sbs B_{c_{n+1},\rho_{n+1}}$. Let $u_n(z):=\kappa_{c_n} (\rho_n (z))$ so that as proved in the first part, $V_{u_n}=d_n$ and $u_n\in C^{\infty}(\ovr{\mathbb{D}})$.

Let
\[u(z):=\sum_{n=0}^{\infty}u_n(z).\] Since $|u_n|\leq |c_n|\leq 2^{-n}$ for all $n$, the sum converges uniformly on $\ovr{\mathbb{D}}$. This shows that $u\in\E$ and properties in a. and c. are satisfied. Using (\ref{Eq:DefnV}) in Theorem \ref{T:CharacHardyClassesBoundaryVal},
\[V_u(\zeta)=\int_{\mathbb{D}}P(z,\zeta)\Delta u(z)=\sum_{n=0}^{\infty}\int_{\mathbb{D}}P(z,\zeta)\Delta u_n(z)=\sum_{n=0}^{\infty}d_n(\zeta)=\psi (\zeta).\]
Due to an equality in Remark \ref{Rem:NormalDer},
\[\int_{\mathbb{D}}\Delta u(z)=\int_{\partial\mathbb{D}}V_ud\nu=\int_{\partial\mathbb{D}}\psi d\nu.\]
Hence $u\in\E_0$ if and only if $\psi\in L^1(d\nu)$. The proof is completed.
\end{proof}

We have now the following converse to Theorem \ref{T:InnerOuterFactorization} to answer the first question in the introduction.
\bT\label{Th:FindingExhaustion} Let $\varphi$ be a zero free analytic function on  $\mathbb{D}$ so that $|\varphi^*|$ equals $\nu$-almost everywhere to an upper semicontinuous function on $\partial\mathbb{D}$. Then
there exists a $u\in\E(\mathbb{D})$
so that $H^p_u=\M_{\varphi, p}$  and we have isometric isomorphism of two spaces.
\eT

\begin{proof}
It is enough to prove the theorem when $p=2$. Since  $|\varphi^*|$ is upper semicontinuous on
$\partial\mathbb{D}$, there exists a constant $m$ so that $|\varphi^*|\leq m$. Hence the function $\psi:=1/|\varphi^*|^2$ is lower semicontinuous and $\psi \geq 1/m$. Let
$u\in\mathcal{E}(\mathbb{D})$ be the exhaustion provided by
Theorem \ref{Th:ConstructExhaustionForH1} for the function
$\psi$ so that $V_u=1/|\varphi^*|^2$.
If $f\in H^2_u$, we write $f=\varphi f_0$, where $f_0=f/\varphi$.
Then
\[\|f\|_{2,u}^2=\int_0^{2\pi}|f(e^{i\theta})|^2 V_{u}(e^{i\theta} )\,d\theta =
\|f\|^2_{\M_{\varphi, 2}}=\|f_0\|_2^2<\infty
.\] Thus $f_0\in H^2$ and we have shown that $H^2_u\sbs \M_{\varphi,
2}$. Conversely, if $f\in  \M_{\varphi, 2}$, then
clearly $f\in H^2_u$ from the same equality above. The
mapping $f\mapsto \varphi  f_0$ is clearly an isomorphism of $H^2_u$
onto $\M_{\varphi, 2}$ which is an isometry.
\end{proof}

When the weight function $V_u$ is smooth enough, there is a connection with the corresponding subharmonic exhaustions and the bi-Laplacian equation $\Delta ^2 u=0$. This is explained in the next result.
\bT\label{Th:ConstructExhaustionForH1} Let $\psi\in C^1(\partial
\mathbb{D})$ be a nonnegative function. Then there exists a function
$u$ and a constant $M$ with the following properties:
\begin{itemize}
\item[a.] $u\in\E_0(\mathbb{D})$ and $u$ is real analytic on $\mathbb{D}$.
\item[b.] $V_u(\zeta)=\frac{\partial u}{\partial n}(\zeta)=\psi (\zeta)+M$ for every $\zeta\in \partial \mathbb{D}$.
\item[c.] $u$ satisfies the bi-Laplacian equation $\Delta^2u=0$ on $\mathbb{D}$.
\end{itemize}
\eT
\begin{proof}
Let $u (z):=\frac{1}{2}(|z|^2-1)[P\psi (z)+M]$, where $P\psi (z)$ is the harmonic extension of $\psi$ on $\mathbb{D}$.
Then using polar coordinates $\Delta u(z)=2[P\psi (z)+M]+2|z|\frac{\partial P\psi (z)}{\partial r}$.
Note that $P\psi\in C^1(\ovr{\mathbb{D}})$. Now take $M$ large enough so that $\Delta u(z)\geq 0$ on $\mathbb{D}$. Again taking the Laplacian it can be checked that $\Delta ^2u(z)=0$. Hence $\Delta u$ is harmonic on $\mathbb{D}$ and since
\[\int_{\mathbb{D}}\Delta u=c(P\psi (0)+M)\leq c\|\psi\|_{\infty}+cM<\infty,\]
$u\in\E_0$. Clearly $u$ is real analytic on $\mathbb{D}$. On $\partial\mathbb{D}$ we have
\[V_u(\zeta)=\frac{\partial u}{\partial r}(\zeta)=\psi (\zeta)+M\]  for every $\zeta\in \partial \mathbb{D}$.
\end{proof}

\begin{Remark}
Equality in (\ref{Eq:OnV}) shows also that
\[2\int_0^{2\pi}\log |z-e^{i\theta}|\left [\psi(e^{i\theta})+\frac{\partial P\psi (e^{i\theta})}{\partial r}+M\right ]d\theta =
\int_{\mathbb{D}}\log |1-\ovr wz|\Delta u(w)\] for every
$z\in\ovr{\mathbb{D}}$. Therefore, in fact, $u$ can be written as
the difference of two potentials
\[u(z)=\frac{1}{2\pi}\int_{\mathbb{D}}\log |z-w|\Delta u(w)dw-\frac{1}{\pi}\int_0^{2\pi}\log |z-e^{i\theta}|
\left [\psi(e^{i\theta})+\frac{\partial P\widetilde{\psi} (e^{i\theta})}{\partial r}+M\right ]d\theta\] for every $z\in\ovr{\mathbb{D}}$. Here $\Delta u$ is harmonic.
\end{Remark}

When $v\in\E (\mathbb{D})$, let $R(v)$ denote the class of all
functions $u\in\E (\mathbb{D})$ which generates the same space
$H^p_v=H^p_u$. We know a "good" representative in $R(v)$ for certain cases as a
consequence of Theorem \ref{Th:ConstructExhaustionForH1}. \bT Let
$v\in\E_0(\mathbb{D})$ so that $V_v$ is bounded and $PV_v+|z|\frac{\partial PV_v}{\partial r}\geq 0$ on $\mathbb{D}$. Then $R(v)$ contains a function $u\in\E_0$ which is
real analytic and satisfies the bi-Laplacian equation $\Delta^2 u=0$
on $\mathbb{D}$. Moreover, $V_u=V_v$ and the weight function $V_u$ can be found by
using the equation
\[V_u(e^{i\theta})=\frac{1}{2}\int_0^1\Delta u(se^{i\theta})ds.\]\eT
\begin{proof}
Let $u(z):=\frac{1}{2}(|z|^2-1)PV_v (z)$. Then $\Delta u(z)=2PV_v (z)+2|z|\frac{\partial PV_v (z)}{\partial r}\geq 0$ by assumption. Hence $u\in \E_0$, $u$ is
real analytic and satisfies the bi-Laplacian equation $\Delta^2 u=0$
on $\mathbb{D}$. Let $h(z):=\Delta u(z)$ and $h_s(z):=h(sz)$ for $0<s<1$. By (\ref{Eq:DefnV}) of Theorem \ref{T:CharacHardyClassesBoundaryVal}
\begin{eqnarray*}V_u(e^{i\theta}) &=&\int_{\mathbb{D}}P(z,e^{i\theta})h(z)dz=\lim_{s\to 1}\int_{\mathbb{D}}P(z,e^{i\theta})h_s(z)dz\\ \nonumber
&=& \lim_{s\to 1}\int_0^1r\int_0^{2\pi}P(re^{it},e^{i\theta})\left [ \frac{1}{2\pi}\int_0^{2\pi}h_s(e^{i\eta})P(re^{it},e^{i\eta})d\eta   \right]dtdr\\ \nonumber
&=& \lim_{s\to 1}\int_0^1r\int_0^{2\pi}h_s(e^{i\eta})\left [ \frac{1}{2\pi}\int_0^{2\pi}P(re^{it},e^{i\eta})P(re^{i\theta},e^{it}) dt  \right]d\eta dr\\ \nonumber
&=& \lim_{s\to 1}\int_0^1r\int_0^{2\pi}h_s(e^{i\eta})P(r^2e^{i\theta},e^{i\eta})d\eta dr\\ \nonumber
&=& \lim_{s\to 1}\int_0^1rh_s(r^2e^{i\theta})dr=\frac{1}{2}\int_0^1\Delta u(re^{i\theta})dr.
\end{eqnarray*}

\end{proof}

\section{Representation of linear functionals}
First we describe the the space of annihilators of $H^p_u$ in
$L^q(Vd\theta)$, where $1/p+1/q=1$. Let $G\in L^q(d\theta)$ and
$f=\varphi^{2/p}F\in H^p_u$. Define
\[L_G(f)=L_G(\varphi^{2/p}F):=\int_0^{2\pi}
F(e^{i\theta})G(e^{i\theta})d\theta.\] Then $L_G$ belongs to
$(H^p_u)^*$ since $|L_G(f)|\leq
\|F\|_p\|G\|_q=\|f\|_{u,p}\|G\|_{q}$. We denote by $H^q_{u,0}$ the
class of functions $g$ in $H^q_u$ with $g(0)=0$. Then $H^q_{u,0}$ is
isometrically isomorphic to $H^q_0$ which is the space of functions
$g\in H^q$ with $g(0)=0$. \bT\label{T:AnnihilatorHpu} For $1\leq
p<\infty$, $(H^p_u)^{\perp}$ is isometrically isomorphic to
$H^q_{0}$ which is isometrically isomorphic to $H^q_{u,0}$ or
$H^q_{u}$.\eT

\begin{proof}
Suppose $g\in L^q(Vd\theta)$ is an annihilator of $H^p_u$. Then
\[\int_0^{2\pi}\varphi^{2/p}(e^{i\theta})g(e^{i\theta})V(e^{i\theta})e^{in\theta}d\theta
=0\] for every $n=0,1,2,\ldots$. Therefore $\varphi^{2/p}gV$ is the
boundary function of some $G\in H^1$ with $G(0)=0$. In fact $G$ is
determined uniquely by $g$. From the equality
\[|g|^qV=|\varphi|^2|G|^qV=|G|^q\] we see that $G\in H^q$ and $\|g\|_{u,q}=\|G\|_q$. Take any $f=\varphi^{2/p}F\in H^p_u$. Then
\[\int_0^{2\pi}f(e^{i\theta})g(e^{i\theta})V(e^{i\theta})d\theta=
\int_0^{2\pi}F(e^{i\theta})G(e^{i\theta})d\theta=0.\]

Conversely, take any $G\in H^q$. Now from \cite[Sec. 7.2]{Duren} if
$G\in H^q_0$, then $L_G\in (H^p_u)^{\perp}$. Hence the map $G\mapsto
L_G$ from $H^q_0$ onto $(H^p_u)^{\perp}$ is an isometric
isomorphism.
\end{proof}
Theorem \ref{T:AnnihilatorHpu} gives a canonical representation of
$(H^p_u)^*$ as in the next statement which can be compared to the
classical case (see \cite[Theorem 7.3]{Duren} for example). \bT For
$1\leq p<\infty$, $(H^p_u)^*$ is isometrically isomorphic to
$L^q(Vd\theta)/H^q_u$. Furthermore, if $1<p<\infty$, for each $L\in
(H^p_u)^*$ there exists a unique $G\in H^q_u$ so that $L(f)=L_G(f)$
for every $f\in H^q_u$. For each $L\in (H^1_u)^*$ there exists a
function $G\in H^{\infty}_u$ so that $L(f)=L_G(f)$ for every $f\in
H^1_u$. \eT

The next theorem describes the preduals of $H^p_u$.

\bT Let $u$ be a subharmonic exhaustion function on $\mathbb{D}$. If
$1 < p \leq \infty$ and $1/p + 1/q = 1$, then:
\begin{itemize}
\item[i.] $H^p_u=\left (L^q_u/H^q_{u,0}\right)^*$.

\item[ii.] $H^p_{u,0}=\left (L^q_u/H^q_u\right)^*$.
\end{itemize}
\eT

\begin{proof}
Let $\Gamma$ be a bounded linear functional on $L^q_u/H^q_{u,0}$.
Then, by composing $\Gamma$ with the canonical projection of $L^q_u$
onto $L^q_u/H^q_{u,0}$, $\Gamma$ gives a linear functional on
$L^q_u$ with the same norm as on $L^q_u/H^q_{u,0}$. So, for any $f
\in L^q_u$, using equation (\ref{Eq:VSign}), Theorem
\ref{Th:IsometryLpuLp} and Corollary \ref{Cor:IsometryHpuHp} we have
\begin{eqnarray}\label{Eq:DefnLinFuncLqQuotientHqu}
\Gamma (f)=\Gamma (f+H^q_{u,0}) = \int_0^{2\pi}
[f(e^{i\theta})\varphi^{-2/q}(e^{i\theta})]
[G(e^{i\theta})\varphi^{-2/p}(e^{i\theta})]d\theta,
\end{eqnarray}
where $G\in L^p_u$ with $\|G\|_{p,u}=\|\Gamma\|$. We have $\Gamma
(e^{in\theta}\varphi^{2/q})=0$ for every integer $n\geq 1$. Hence
$G\in H^p_u$. Conversely, any function $G\in H^p_u$ gives rise to a
linear functional $\Gamma$ on $L^q_u/H^q_{u,0}$ by formula
(\ref{Eq:DefnLinFuncLqQuotientHqu}). This proves the first
assertion. The second part is proved by a similar argument.
\end{proof}

\section{Extremal problems}
We are now ready to discuss the related extremal problems. For fixed
$g\in L^q(Vd\theta)$ the extremal problem is to find
\begin{eqnarray}\label{Eq:SupProblemHpu}
\Lambda (g):= \sup\left\{|\lambda (f)|:f\in H^p_u,\, \|f\|_{p,u}\leq
1\right\},
\end{eqnarray}
where
\begin{eqnarray}
\lambda (f):= \frac{1}{2\pi i}\int_{|z|=1} F(z)G(z)dz =
\frac{1}{2\pi }\int_{0}^{2\pi}
f(e^{i\theta})g(e^{i\theta})V(e^{i\theta})e^{i\theta}d\theta
\end{eqnarray}
and we use the correspondence $f = \varphi^{2/p}F$, $g =
\varphi^{2/q}sgn(\varphi^2)G$ provided by Theorem
\ref{Th:IsometryLpuLp} and Corollary \ref{Cor:IsometryHpuHp}. The
related dual extremal problem is to find the function $g_0\in H^q_u$
so that
\begin{eqnarray}\label{Eq:InfProblemHqu}
\Gamma(g):=\inf \left\{\|g-h\|_{q,u}:h\in
H^q_u\right\}=\|g-g_0\|_{q,u}.
\end{eqnarray}

The proof of the following existence and uniqueness theorem for the
extremal problems follows in view of Theorem \ref{Th:IsometryLpuLp},
Corollary \ref{Cor:IsometryHpuHp} and \cite[Theorem 8.1]{Duren}.

\bT Let $1\leq p\leq \infty$, $1/p+1/q=1$ and $g\in L^q(Vd\theta)$.
\begin{itemize}
\item[i.] The duality relation $\Lambda (g)=\Gamma(g)$ holds.

\item[ii.] If $p>1$, there is a unique extremal function $f\in
H^p_u$ for which $\lambda (f)>0$. The dual extremal problem has a
unique solution.

\item[iii.] If $p=1$ and $G(e^{i\theta})$ is continuous, at least
one solution to the extremal problem exists. If $p=1$, the dual
extremal problem has at least one solution; it is unique if the
extremal problem has a solution.
\end{itemize}
\eT


\bibliographystyle{amsplain}

\end{document}